\documentstyle[12pt, amsfonts]{amsart}

\def\R{{\mathbb R}}

\def\C{{\mathbb C}}

\newtheorem{df}{Definition}
\newtheorem{theorem}{Theorem}[section]
\newtheorem{lemma}[theorem]{Lemma}
\newtheorem{pr}[theorem]{Proposition}
\newtheorem{co}[theorem]{Corollary}

\def\N{{\mathbb N}}

\def\vol{\mbox{\rm Vol}}

\def\e{\varepsilon}

\def\pf {{\sl Proof. }}
\def\endpf{ \begin{flushright}
$ \Box $ \\
\end{flushright}}
\begin{document}
 
\title[Complex intersection bodies]{Complex intersection bodies}

\author{A.\ Koldobsky} 
\address{Department of Mathematics\\ University of Missouri\\ Columbia\\ Missouri 65211}
\email{koldobskiya\@ missouri.edu}

\author{G.\ Paouris} 
\address{Department of Mathematics\\Texas A $\&$ M University\\College Station \\TX 77843- 3368\\}
\email{grigoris\@ math.tamu.edu}

\author{M.\ Zymonopoulou}\address{ Department of Mathematics\\University of Crete\\ Heraklio\\Crete}
\email{marisa.zym\@ gmail.com}
 
\begin{abstract} We introduce complex intersection bodies and show that
their properties and applications are similar to those of their real counterparts.
In particular, we generalize Busemann's theorem to the complex case  by proving  
that complex intersection bodies of symmetric complex convex bodies are also
convex. Other results include stability in the complex Busemann-Petty problem for arbitrary 
measures and the corresponding hyperplane inequality for measures of complex 
intersection bodies.
\end{abstract}
\maketitle

\section{Introduction}

The concept of an intersection body was introduced by Lutwak \cite{L}, as
part of his dual Brunn-Minkowski theory. In particular, these bodies 
played an important role in the solution of the Busemann-Petty problem. 
Many results on intersection bodies have appeared in recent years (see \cite{G,K1,KY}
and references there), but almost all of them apply to the real case. The goal of this paper 
is to extend the concept of an intersection body to the complex case. 

Let $K$ and $L$ be origin symmetric star bodies in $\R^n.$ Following \cite{L}, 
we say that $K$ is the {\it intersection body of} $L$ if the radius of $K$ in every direction
is equal to the volume of the central hyperplane section of $L$ perpendicular
to this direction, i.e. for every $\xi\in S^{n-1},$
\begin{equation}\label{intbodyofstar}
\|\xi\|_K^{-1} = |L\cap \xi^\bot|,
\end{equation}
where $\|x\|_K=min\{a\ge 0:\ x\in aK\}$, $\xi^\bot=\{x\in \R^n:\ (x,\xi)=0\},$
and $|\cdot|$ stands for volume.
By a theorem of Busemann \cite{Bu1} the intersection body of an origin symmetric
convex body is also convex. However, intersection bodies of convex bodies 
form just a small part of the class of intersection bodies. In particular,  by results of  Hensley 
\cite{H} and Borell \cite{Bor}, intersection bodies of symmetric convex bodies are isomorphic 
to an ellipsoid,  i.e.  $d_{BM}(I(K), B_{2}^{n})\le c$ where $d_{BM}$ is the Banach-Mazur 
distance and $c>0$ is a universal constant. 

The right-hand
side of (\ref{intbodyofstar}) can be written using the polar formula for volume:
$$\|\xi\|_K^{-1}= \frac{1}{n-1} \int_{S^{n-1}\cap \xi^\bot} \|\theta\|_L^{-n+1} d\theta=
\frac{1}{n-1} {\cal{R}}(\|\cdot\|_L^{-n+1})(\xi),$$
where the operator ${\cal{R}}: C(S^{n-1})\to C(S^{n-1})$ is the spherical Radon transform defined by
$${\cal{R}}f(\xi) = \int_{S^{n-1}\cap \xi^\bot} f(x) dx.$$ 
This means that a star body $K$ is the intersection body of a star body if and only if the function $\|\cdot\|_K^{-1}$
is the spherical Radon transform of a continuous positive function on $S^{n-1}.$

A more general class of bodies was introduced in \cite{GLW}. A star body
$K$ in $\R^n$ is called an {\it intersection body} if there exists a finite Borel measure \index{intersection body}
$\mu$ on the sphere $S^{n-1}$ so that $\|\cdot\|_K^{-1}= R\mu$ as functionals on 
$C(S^{n-1}),$ i.e. for every continuous function $f$ on $S^{n-1},$
\begin{equation} \label{defintbody}
\int_{S^{n-1}} \|x\|_K^{-1} f(x)\ dx = \left( {\cal{R}}\mu, f \right)= \int_{S^{n-1}} Rf(x)\ d\mu(x).
\end{equation}

We introduce complex intersection bodies along the same lines. In Section \ref{CIBofstarbodies}
we define complex intersection bodies of complex star bodies, and in  Section \ref{bus} we study
complex intersection bodies of convex bodies. While the complex version of Busemann's theorem 
requires a serious effort, the extension of the Hensley-Borell theorem to the complex case
follows from a result of Ball \cite{Ball1}. In Section \ref{radonandfourier} we prove 
that the complex spherical Radon transform and the Fourier transform of distributions coincide (up to a constant)
on a class of $(-2n+2)$-homogeneous functions on $\R^{2n}$ with symmetries determined by
the complex structure. A similar result in the real case was crucial for the study of real intersection 
bodies. We use this result in Section \ref{CIB}, where we define complex intersection bodies and prove 
a Fourier characterization of intersection bodies: an origin symmetric  complex star body $K$ in $\R^{2n}$
is a complex intersection body if and only if the function $\|\cdot\|_K^{-2}$ represents a positive definite
distribution. We use this characterization in Section \ref{Char of CIB} to show that the class of complex intersection bodies 
coincides with the class of real 2-intersection bodies in $\R^{2n}$ and, at the same time, with the class
of generalized 2-intersection bodies, provided that bodies from the real classes possess symmetries
determined by the complex structure of $\R^{2n}.$ 
The latter allows to extend to the complex case a result of Goodey
and Weil \cite{GW} by showing that all symmetric complex intersection
bodies can be obtained as limits in the radial metric of complex radial sums of ellipsoids. Finally, Section
\ref{stab-hyper} deals with stability in the complex Busemann-Petty problem for arbitrary
measures and related hyperplane inequalities.


\section{Complex intersection bodies of star bodies}\label{CIBofstarbodies}

The theory of real convex bodies goes back to ancient times and 
continues to be a very active field now. The situation with complex convex bodies 
is different, as no systematic studies of these bodies have been carried out, and results appear
only occasionally; see  for example \cite{KKZ, KZ, AB, Ru1, Zy, Zy1}.

origin symmetric convex bodies in $\C^n$ are the unit balls of norms on $\C^n.$
We denote by $\|\cdot\|_K$
the norm corresponding to the body $K:$
$$K=\{z\in \C^n:\ \|z\|_K\le 1\}.$$
In order to define volume, we identify $\C^n$ with $\R^{2n}$ using the standard mapping
$$\xi = (\xi_1,...,\xi_n)=(\xi_{11}+i\xi_{12},...,\xi_{n1}+i\xi_{n2})
 \mapsto  (\xi_{11},\xi_{12},...,\xi_{n1},\xi_{n2}).$$
 Since norms on $\C^n$ satisfy the equality
$$\|\lambda z\| = |\lambda|\|z\|,\quad \forall z\in \C^n,\  \forall\lambda \in \C,$$
origin symmetric complex convex bodies correspond to those origin symmetric convex bodies
$K$  in $\R^{2n}$ that are invariant
 with respect to any coordinate-wise two-dimensional rotation, namely for each $\theta\in [0,2\pi]$
 and each $\xi= (\xi_{11},\xi_{12},...,\xi_{n1},\xi_{n2})\in \R^{2n}$
  \begin{equation} \label{rotation}
  \|\xi\|_K =
 \|R_\theta(\xi_{11},\xi_{12}),...,R_\theta(\xi_{n1},\xi_{n2})\|_K,
 \end{equation}
 where $R_\theta$ stands for  the counterclockwise rotation of $\R^2$ by the angle
 $\theta$ with respect to the origin. We shall say that $K$ is a {\it complex convex body
 in $\R^{2n}$} if $K$ is a convex body and satisfies equations (\ref{rotation}).
 
 A compact set $K$ in $\R^n$ is called a star body if the origin is an interior point of $K$, every straight
line passing through the origin crosses the boundary of $K$ at exactly two points and the Minkowski
functional of $K$ defined by
$$\|x\|_K= \min\{a\ge 0:\ x\in aK\},\qquad \forall x\in \R^n$$
is a continuous function on $\R^n.$ The radial function of $K$ is given by 
$$\rho_K(x) = \max\{a>0:\ ax\in K\}.$$
If $x\in S^{n-1},$ then $\rho_K(x)$ is the radius of $K$ in the direction of $x.$
Note that for any unit vector $\xi$, $\rho_K(\xi)=\|\xi\|_K^{-1}.$ The radial metric
in the class of star bodies is defined by
$$\rho(K,L)=\max_{\xi\in S^{n-1}} |\rho_K(\xi)-\rho_L(\xi)|.$$
If the Minkowski functional of a star body $K$ in $\R^{2n}$ is $R_\theta$-invariant 
(i.e. satisfies equations (\ref{rotation})), we say that $K$ is  a {\it complex star body}
in $\R^{2n}.$
 
 For $\xi\in \C^n,
|\xi|=1,$ denote by
$$H_\xi = \{ z\in \C^n:\ (z,\xi)=\sum_{k=1}^n z_k\overline{\xi_k} =0\}$$
the complex hyperplane through the origin, perpendicular to $\xi.$
 Under the standard mapping from $\C^n$ to $\R^{2n}$ the hyperplane $H_\xi$ 
 turns into a $(2n-2)$-dimensional subspace of $\R^{2n}$ orthogonal to the vectors
 $$\xi=(\xi_{11},\xi_{12},...,\xi_{n1},\xi_{n2})\quad  {\rm and } \quad
 \xi^\perp =(-\xi_{12},\xi_{11},...,-\xi_{n2},\xi_{n1}).$$
 The orthogonal two-dimensional subspace $H_\xi^\bot$ has orthonormal basis $\left\{\xi,\xi^\bot\right\}.$
 A star (convex) body $K$ in $\R^{2n}$ is a complex star (convex) body if and only if,
for every $\xi\in S^{2n-1},$  the section $K\cap H_\xi^\bot$ is a two-dimensional Euclidean circle
with radius $\rho_K(\xi)= \|\xi\|_K^{-1}.$

 We introduce complex intersection bodies of complex star bodies using a definition 
 under which these bodies play the same role in complex convexity, as their real counterparts 
 in the real case. We use the notation $|K|$ for the volume of $K;$ the dimension where we consider 
 volume is clear in every particular case. 
 
 \begin{df} \label{intof} Let  $K, L$ be origin symmetric complex star bodies in $\R^{2n}.$ We say that $K$ 
 is the complex intersection 
 body of $L$ and write $K=I_c(L)$ if for every $\xi \in \R^{2n}$
 \begin{equation} \label{def-intbody}
 |K\cap H_\xi^{\perp}|= |L\cap H_\xi|.
 \end{equation}
  \end{df}
  Since $K\cap H_\xi^{\perp}$ is the two-dimensional Euclidean circle with radius $\|\xi\|_K^{-1}$,
  (\ref{def-intbody}) can be written as
  \begin{equation}\label{intbody}
 \pi  \|\xi\|_{I_c(L)}^{-2} =  |L\cap H_\xi|.
  \end{equation}
 All the bodies $K$ that appear as complex intersection bodies of different complex star bodies
form {\it the class of complex intersection bodies of star bodies}. In Section \ref{CIB}, we will introduce
a more general class of {\it complex intersection bodies}.


   \section{The Radon and Fourier transforms of $R_\theta$-invariant functions} \label{radonandfourier}

  Denote by $C_c(S^{2n-1})$ the space of $R_\theta$-invariant continuous functions, i.e.
  continuous real-valued functions $f$ on the unit sphere $S^{2n-1}$ in $\R^{2n}$ satisfying 
  $f(\xi)=f(R_\theta(\xi))$ for all $\xi\in S^{2n-1}$ and all $\theta\in [0,2\pi].$ The complex spherical
  Radon transform is an operator ${\cal{R}}_c: C_c(S^{2n-1})\to C_c(S^{2n-1})$ defined by
  $${\cal{R}}_cf(\xi) = \int_{S^{2n-1}\cap H_\xi} f(x) dx.$$
  Writing volume in polar coordinates, we get that for every complex star body $L$ in $\R^{2n}$
  and every $\xi\in S^{2n-1},$
  \begin{equation}\label{v-r}
  |L\cap H_\xi| = \frac1{2n-2} \int_{S^{2n-1}\cap H_\xi} \|x\|_L^{-2n+2} dx 
  =  \frac1{2n-2} {\cal{R}}_c\left(\|\cdot\|_L^{-2n+2}\right)(\xi),
  \end{equation}
  so the condition (\ref{intbody}) reads as
  \begin{equation}\label{int-radon}
  \|\xi\|_{I_c(L)}^{-2} =  \frac1{2\pi(n-1)} {\cal{R}}_c\left(\|\cdot\|_L^{-2n+2}\right)(\xi).
  \end{equation}
  This means that a complex star body $K$ is a complex intersection body of a star body if and only
  if the function $\|\cdot\|_K^{-2}$ is the complex spherical Radon transform of a continuous positive
  $R_\theta$-invariant function on $S^{2n-1}.$ We use this observation in Section \ref{CIB}, where  we introduce 
  a more general class of complex intersection bodies (not depending on the underlying star body),
  like it was done in the real case in \cite{GLW}. But before that we need several facts connecting the Radon transform
  to the Fourier transform in the complex setting.
  
 We use the techniques of the Fourier transform of distributions; see \cite{GS} for details.
As usual, we denote by ${\cal{S}}(\R^n)$
the Schwartz space of rapidly decreasing infinitely differentiable
functions (test functions) in $\R^n,$ and
${\cal{S}}^{'}(\R^n)$ is
the space of distributions over ${\cal{S}}(\R^n).$

Suppose that $f$ is a locally 
integrable complex-valued function on $\R^n$ with {\it power growth at infinity}, 
i.e. there exists a number $ \beta>0$ so that 
$$\lim_{|x|_2\to \infty} \frac{f(x)} {|x|_2^\beta}=0,$$
where $|\cdot|_2$ stands for the Euclidean norm on $\R^n.$ 
Then $f$ represents a distribution acting by integration: 
for every $\phi\in {\mathcal S},$ 
$$\langle f, \phi \rangle = \int_{\R^n} f(x) \phi(x)\ dx.$$

The Fourier transform of a
distribution $f$ is defined by $\langle\hat{f}, \phi\rangle= \langle f, \hat{\phi} \rangle$ for
every test function $\phi.$ If $\phi$ is an even test function, then $(\hat{\phi})^\wedge = (2\pi)^n \phi$,
so the Fourier transform is self-invertible (up to a constant) for even distributions.

A distribution $f$ is called even homogeneous of degree $p\in \R$  if
$$\langle f(x), \phi(x/\alpha) \rangle = |\alpha|^{n+p} 
\langle f,\phi \rangle$$
for every test function $\phi$ and every 
$\alpha\in \R,\ \alpha\neq 0.$  The Fourier transform of an even
homogeneous distribution of degree $p$ is an even homogeneous
distribution of degree $-n-p.$ 

We say that a distribution is  {\it positive definite} if its Fourier transform is a positive distribution in
the sense that $\langle \hat{f},\phi \rangle \ge 0$ for every non-negative test function $\phi.$
Schwartz's generalization of Bochner's theorem (see, for example, 
\cite[p.152]{GV}) states that a distribution is positive definite
if and only if it is the Fourier transform of 
a tempered measure on $\R^n$. Recall that \index{tempered measure}
a (non-negative, not necessarily finite) measure $\mu$ is  
called tempered if 
$$\int_{\R^n} (1+|x|_2)^{-\beta}\ d\mu(x)< \infty$$
for some $\beta >0.$ 

Our definition of a star body $K$
assumes that the origin is an interior point of $K.$ If  $0<p<n,$
then $\|\cdot\|_K^{-p}$  is a locally integrable function on $\R^n$ and represents an even 
homogeneous of degree $-p$ distribution. If $\|\cdot\|_K^{-p}$ represents a positive definite
distribution for some $p\in (0,n),$ then its Fourier transform is a tempered measure which 
is at the same time a homogeneous distribution of degree $-n+p.$ One can express such 
a measure in polar coordinates, as follows.

\begin{pr} \label{posdef}  (\cite[Corollary 2.26]{K1}) 
Let $K$ be an origin symmetric convex body in $\R^n$ and
$p\in (0,n).$ The function $\|\cdot\|_K^{-p}$ represents a 
positive definite distribution on $\R^n$ if and only if there 
exists a finite Borel measure $\mu$ on $S^{n-1}$ so that 
for every even test function $\phi,$
$$
\int_{\R^n} \|x\|_K^{-p} \phi(x)\ dx = \int_{S^{n-1}} \left(
\int_0^\infty t^{p-1} \hat\phi(t\xi) dt \right) d\mu(\xi).
$$
\end{pr}

For any even continuous function $f$ on the sphere $S^{n-1}$ and any non-zero number $p\in \R,$ 
we denote by $f\cdot r^{p}$ the  extension of $f$ to an even homogeneous function of degree $p$ 
on $\R^n$ defined as follows. If $x\in \R^n,$ then  $x=r\theta,$ where $r=|x|_2$ and $\theta= x/|x|_2.$ We put
$$f\cdot r^p(x) = f\left(\theta\right) r^p.$$
It was proved in \cite[Lemma 3.16]{K1} that, for any $p\in (-n,0)$ and infinitely smooth function $f$
on $S^{n-1},$ the Fourier transform of $f\cdot r^{-p}$ 
is equal to another infinitely smooth function $h$ on $S^{n-1}$ extended to an even homogeneous 
of degree $-n+p$ function $h\cdot r^{-n+p}$ on the whole of $\R^n.$
The following Parseval formula on the sphere was proved in \cite[Corollary 3.22]{K1}.
\begin{pr} Let $f,g$ be even infinitely smooth functions on $S^{n-1},$
and $p\in (0,n).$ Then
\begin{equation} \label{parseval}
\int_{S^{n-1}} (f\cdot r^{-p})^\wedge(\theta) (g\cdot r^{-n+p})^\wedge(\theta) =
(2\pi)^n \int_{S^{n-1}}  f(\theta)g(\theta) \ d\theta. 
\end{equation}
\end{pr}

We need a simple observation that will, however, provide the basis for applications of the 
Fourier transform to complex bodies.
\begin{lemma} \label{const} Suppose that $f\in C_c(S^{2n-1})$ is an even infinitely smooth
function. Then for every $0<p<2n$ and $\xi\in S^{2n-1}$ the Fourier transform of the distribution $f\cdot r^{-p}$ is
a constant function on $S^{2n-1}\cap H_\xi^\bot.$
\end{lemma}
\pf  By \cite[Lemma 3.16]{K1}, the Fourier transform of $f\cdot r^{-p}$ is a continuous
function outside of the origin in $\R^{2n}.$ The function $f$ is invariant with respect
to all $R_\theta$, so by the
connection between the Fourier transform of distributions and linear transformations,
the Fourier transform of $f\cdot r^{-p}$ is also invariant with respect to all $R_\theta.$
Recall that the two-dimensional space $H_\xi^\bot$ is spanned by  vectors
$\xi$ and $\xi^\bot$ (see the Introduction). Every vector in $S^{2n-1}\cap H_\xi^\bot$
is the image of $\xi$ under one of the coordinate-wise rotations $R_\theta$, so
the Fourier transform of $f\cdot r^{-p}$ is a constant function on $S^{2n-1}\cap H_\xi^\bot.$
\endpf

The following connection between the Fourier and Radon transforms is well-known; see 
for example \cite[Lemma 3.24]{K1}.
\begin{pr} \label{f-r} Let $1\le k <n,$ and let $\phi\in {\cal{S}}(\R^n)$ be  an even test function.
Then for any $(n-k)$-dimensional subspace $H$ of $\R^n$
$$\int_H \phi(x) dx = \frac1{(2\pi)^k} \int_{H^\bot} \hat{\phi}(x) dx.$$
\end{pr}
We also use the spherical version of Proposition \ref{f-r}; see \cite[Lemma 3.25]{K1}. 
\begin{pr} \label{perp} Let $\phi$ be an even infinitely smooth function on $S^{n-1},$
let $0<k<n,$ and let $H$ be an arbitrary $(n-k)$-dimensional subspace of $\R^n.$
Then
$$\int_{S^{n-1} \cap H} \phi(\theta)d\theta =
 \frac1{(2\pi)^k} \int_{S^{n-1}\cap H^\bot} \left(\phi\cdot r^{-n+k}\right)^\wedge(\theta) d\theta.$$
\end{pr}
Let us translate the latter fact to the complex situation.
\begin{lemma} \label{radon-fourier} Let $\phi\in C_c(S^{2n-1})$ be an even infinitely smooth function. Then for every
$\xi\in S^{2n-1}$
$$ {\cal{R}}_c\phi(\xi)= \frac1{2\pi} \left(\phi\cdot r^{-2n+2}\right)^\wedge (\xi).$$
\end{lemma}
\pf By Proposition \ref{perp},
$${\cal{R}}_c\phi(\xi) = \int_{S^{2n-1} \cap H_\xi} \phi(\theta)d\theta= \frac1{(2\pi)^2} \int_{S^{2n-1}\cap H_\xi^\bot}
 \left(\phi\cdot r^{-2n+2}\right)^\wedge (\theta) d\theta,$$
 and, by Lemma \ref{const}, the function under the integral in the right-hand side is constant
 on $S^{2n-1}\cap H_\xi^\bot.$ The value of this constant is the function value at 
 $\xi \in S^{2n-1}\cap H_\xi^\bot.$ Also, recall that $S^{2n-1}\cap H_\xi^\bot$ is the two-dimensional
 Euclidean unit circle, so
 $$\int_{S^{2n-1}\cap H_\xi^\bot} \left(\phi\cdot r^{-2n+2}\right)^\wedge (\theta) d\theta =
 2\pi  \left(\phi\cdot r^{-2n+2}\right)^\wedge (\xi). $$
\endpf

\begin{lemma} \label{self-dual}The complex spherical Radon transform is self-dual, i.e. for any even functions
$f,g\in C_c(S^{2n-1}),$
$$\int_{S^{2n-1}} {\cal{R}}_cf(\xi) g(\xi) d\xi = \int_{S^{2n-1}} f(\theta) {\cal{R}}_c g(\theta) d\theta.$$ 
\end{lemma}

\pf By approximation, it is enough to consider the case where $f,g$ are infinitely smooth. 
For some infinitely smooth even function $h\in C_c(S^{2n-1}),$ we have $g\cdot r^{-2n+2}= \left(h\cdot r^{-2}\right)^\wedge,$
then $ \left(g\cdot r^{-2n+2}\right)^\wedge= (2\pi)^{2n} h\cdot r^{-2}.$
By Lemma \ref{radon-fourier} and the spherical Parseval formula  (\ref{parseval}), 
$$\int_{S^{2n-1}} {\cal{R}}_cf(\xi) g(\xi) d\xi = \frac1{2\pi} \int_{S^{2n-1}} \left(f\cdot r^{-2n+2}\right)^\wedge(\xi)
(g\cdot r^{-2n+2}) (\xi) d\xi$$
$$= \frac{(2\pi)^{2n}}{2\pi} \int_{S^{2n-1}} \left(f\cdot r^{-2n+2}\right)^\wedge(\xi)
(h\cdot r^{-2})^\wedge (\xi) d\xi$$
$$= \frac1{2\pi} \int_{S^{2n-1}} f(\theta) \left(g\cdot r^{-2n+2}\right)^\wedge(\theta) d\theta = 
\int_{S^{2n-1}} f(\theta) {\cal{R}}_c g(\theta) d\theta.$$
\endpf

We now prove Lemma \ref{radon-fourier} without smoothness assumption. This
result is a complex version of \cite[Lemma 3.7]{K1}.
We say that a distribution $f$ on $\R^{2n}$ is $R_\theta$-invariant  if 
$\langle f, \phi(R_\theta(\cdot))\rangle = \langle f, \phi \rangle$ for every
test function $\phi \in {\cal{S}}(\R^{2n})$ and every $\theta\in [0,2\pi].$
If $f$ and $g$ are $R_\theta$-invariant distributions, and 
$\langle f, \phi \rangle=\langle g, \phi \rangle$ for any test function $\phi$
that is invariant with respect to all $R_\theta,$ then $f=g.$ This follows from
the observation that the value of an $R_\theta$-invariant distribution on a test function $\phi$
does not change if $\phi$ is replaced by the function $\frac1{2\pi}\int_0^{2\pi} \phi(R_\theta(\cdot))d\theta.$

\begin{lemma} \label{r-s} Let $f\in C_c(S^{2n-1})$ be an even function. Then the Fourier transform of $f\cdot r^{-2n+2}$
is a continuous function on the sphere extended to a homogeneous function of degree -2 on the whole $\R^{2n}.$
Moreover, on the sphere this function is equal (up to a constant) to the complex spherical Radon transform of $f$:
for any $\xi\in S^{2n-1},$
$${\cal{R}}_cf(\xi) = \frac1{2\pi} \left(f\cdot r^{-2n+2}\right)^\wedge(\xi).$$
\end{lemma}
\pf  Let $\phi \in {\cal{S}}(\R^{2n})$ be any even $R_\theta$-invariant test function. Then $\hat{\phi}$ is also
an even $R_\theta$-invariant test function. By Lemma \ref{const}, for any $\xi\in S^{2n-1},$
\begin{equation}\label{eq1}
 \int_{H_\xi^\bot} \hat{\phi}(x) dx = 
\int_{S^{2n-1}\cap H_\xi^\bot} \left(\int_0^\infty r\hat{\phi}(r\theta)dr \right) d\theta =
2\pi \int_0^\infty r\hat{\phi}(r\xi)dr.
\end{equation}
By Proposition \ref{f-r} and Lemma \ref{self-dual},
$$\langle \left(f\cdot r^{-2n+2}\right)^\wedge, \phi \rangle = \int_{\R^{2n}} |x|_2^{-2n+2} f(x/|x|_2) \hat{\phi}(x) dx$$
$$= \int_{S^{2n-1}} f(\xi) \left(\int_0^\infty r\hat{\phi}(r\xi) dr \right) d\xi = 
\frac1{2\pi} \int_{S^{2n-1}} f(\xi) \left(\int_{H_\xi^\bot} \hat{\phi}(x) dx\right) d\xi$$
$$ = 2\pi \int_{S^{2n-1}} f(\xi) \left(\int_{H_\xi} \phi(x) dx \right) d\xi$$ 
$$= 2\pi \int_{S^{2n-1}} f(\xi) \left(\int_{S^{2n-1}\cap H_\xi} \left(\int_0^\infty r^{2n-3} \phi(r\theta) dr \right) d\theta \right) d\xi$$
$$= 2\pi \int_{S^{2n-1}} f(\xi)\  {\cal{R}}_c\left(\int_0^\infty r^{2n-3}\phi(r\cdot)dr\right)(\xi) d\xi$$
$$= 2\pi \int_{S^{2n-1}} \left(\int_0^\infty r^{2n-3}\phi(r\theta)dr\right) {\cal{R}}_cf(\theta) d\theta$$
$$= 2\pi \int_{\R^{2n}} |x|_2^{-2} {\cal{R}}_cf(x/|x|_2) \phi(x) dx.$$
We get that for every even $R_\theta$-invariant test function $\phi,$
$$\langle \left(f\cdot r^{-2n+2}\right)^\wedge, \phi \rangle = 2\pi \langle |x|_2^{-2}  {\cal{R}}_cf(x/|x|_2) , \phi \rangle,$$
so even $R_\theta$-invariant distributions $ \left(f\cdot r^{-2n+2}\right)^\wedge$ and 
$2\pi |x|_2^{-2}  {\cal{R}}_cf(x/|x|_2)$ are equal.
\endpf

Lemma \ref{r-s} implies the following Fourier transform formula for the volume of sections of star bodies.
Note that the real version of this formula was proved in \cite{K2}, and that the complex formula below
was proved in \cite{KKZ} for infinitely smooth bodies by a different method;  here
we remove the smoothness condition.

 \begin{theorem} \label{volume-ft} Let $K$ be an origin symmetric
complex star body in $\R^{2n}, n\ge 2.$  For every $\xi\in S^{2n-1},$ we have
$$|K\cap H_\xi| = \frac{1}{4\pi(n-1)} \left(\|x\|_K^{-2n+2}\right)^\wedge(\xi).$$
\end{theorem}
\pf By (\ref{v-r}) and Lemma \ref{r-s} applied to the function $f(\theta)= \|\theta\|_K^{-2n+2},$
$$|K\cap H_\xi| = \frac1{2n-2} {\cal{R}}_c\left(\|\cdot\|_K^{-2n+2}\right)(\xi)= \frac{1}{4\pi(n-1)} \left(\|x\|_K^{-2n+2}\right)^\wedge(\xi).$$
\endpf

We use Theorem \ref{volume-ft} to prove the complex version of the Minkowski-Funk theorem
saying that an origin symmetric star body is uniquely determined by volume of its central
hyperplane sections; see \cite[Corollary 3.9]{K1}.

\begin{co} If $K,L$ are origin symmetric complex star bodies in $\R^{2n},$ and 
their intersection bodies $I_c(K)$ and $I_c(L)$  coincide, then $K=L.$ 
\end{co}
\pf The equality of intersection bodies means that, for every $\xi\in S^{2n-1},$
$|K\cap H_\xi|= |L\cap H_\xi|.$ By Theorem \ref{volume-ft}, homogeneous of degree -2
continuous on $\R^{2n}\setminus \{0\}$ functions $\left(\|\cdot\|_K^{-2n+2}\right)^\wedge$ and $\left(\|\cdot\|_L^{-2n+2}\right)^\wedge$
coincide on the sphere $S^{2n-1}$, so they are also equal as distributions on the whole of $\R^{2n}.$ The result
follows from the uniqueness theorem for the Fourier transform of distributions.
\endpf


\section{Complex intersection bodies} \label{CIB}

We are going to define the class of complex intersection bodies by extending the equality (\ref{int-radon}) to measures,
as it was done in the real case in \cite{GLW}. We say that a finite Borel measure $\mu$ on
$S^{2n-1}$ is $R_\theta$-invariant if for any continuous function $f$ on $S^{2n-1}$ and any $\theta\in [0,2\pi]$,
$$\int_{S^{2n-1}} f(x) d\mu(x) = \int_{S^{2n-1}} f(R_\theta x) d\mu(x).$$
The complex spherical Radon transform of an $R_\theta$-invariant measure $\mu$ is defined
as a functional ${\cal{R}}_c\mu$ on the space $C_c(S^{2n-1})$ acting by 
$$ \left({\cal{R}}_c\mu, f \right) = \int_{S^{2n-1}} {\cal{R}}_cf(x) d\mu(x).$$

\begin{df}  \label{int} An origin symmetric complex star body $K$ in $\R^{2n}$ is called a {complex intersection body} if there
exists a finite Borel $R_\theta$-invariant measure $\mu$ on $S^{2n-1}$ so that
$\|\cdot\|_K^{-2}$ and ${\cal{R}}_c\mu$ are equal as functionals on $C_c(S^{2n-1}),$ i.e.
for any $f\in C_c(S^{2n-1}),$
$$\int_{S^{2n-1}} \|x||_K^{-2} f (x)\ dx = \int_{S^{2n-1}} {\cal{R}}_c f(\theta) d\mu(\theta).$$
\end{df}

Clearly, ${\cal{R}}_c\mu$ is a finite Borel $R_\theta$-invariant measure on $S^{2n-1}.$
Also, an easy consequence of self-duality of the complex spherical Radon transform (see Lemma \ref{self-dual})
is that if $\mu$ has density $f$ on $S^{2n-1}$, then the measure ${\cal{R}}_c\mu$ has density ${\cal{R}}_cf.$
The latter, in conjunction with (\ref{int-radon}), immediately implies that every complex intersection body of a star body is a complex 
intersection body in the sense of Definition \ref{int}.

Many results on real intersection bodies depend on the following Fourier characterization (see \cite[Theorem 1]{K5}):
an origin symmetric star body $K$ in $\R^n$ is an intersection body if and only if the function $\|\cdot\|_K^{-1}$
represents a positive definite distribution. Complex intersection bodies admit a similar characterization.
To see the connection with the Fourier transform, combine the definition of the intersection body of a star body (\ref{def-intbody}) 
with the result of Theorem \ref{volume-ft}: for every $\xi\in S^{2n-1},$
$$\|\xi\|_{I_c(L)}^{-2} =  \frac1{\pi}|L\cap H_\xi| = \frac{1}{4\pi^2(n-1)} \left(\|x\|_L^{-2n+2}\right)^\wedge(\xi).$$
Both sides of the latter equality are even homogeneous functions of degree -2, so these functions are equal as distributions 
on the whole of $\R^{2n}.$ Since the Fourier transform of even distributions is self-invertible (up to a constant), we get
$$\left(\|\cdot\|_{I_c(L)}^{-2}\right)^\wedge = \frac{(2\pi)^{2n}}{4\pi^2(n-1)} \|\cdot\|_L^{-2n+2} >0,$$
so the distribution $\|\cdot\|_{I_c(L)}$ is positive definite. Moreover, if the Fourier transform  of $\|\cdot\|_K^{-2}$
is an even strictly positive $R_\theta$-invariant function on the sphere, then one can use the latter equality 
to construct a complex star body $L$ such that $K=I_c(L).$ This connection holds for arbitrary 
complex intersection bodies, as shown in the following theorem.

\begin{theorem} \label{posdef} An origin symmetric complex star body $K$ in $\R^{2n}$ is a complex 
intersection body  if and only if the function $\|x\|_K^{-2}$ represents a positive definite distribution on $\R^{2n}.$
\end{theorem}
\pf Suppose that $K$ is a complex intersection body with the corresponding measure $\mu.$ 
To prove that $\|\cdot\|_K^{-2}$ is a positive definite distribution,
it is enough to show that $\langle (\|\cdot\|_K^{-2})^\wedge, \phi\rangle \ge 0$ for every even $R_\theta$-invariant
non-negative test function $\phi.$ By Definition \ref{int} and Proposition \ref{f-r},
$$\langle (\|\cdot\|_K^{-2})^\wedge, \phi\rangle = \int_{\R^{2n}} \|x\|_K^{-2}\hat{\phi}(x) dx$$
$$= \int_{S^{2n-1}} \|\theta\|_K^{-2} \left(\int_0^\infty r^{2n-3} \hat{\phi}(r\theta) dr \right) d\theta$$
$$= \int_{S^{2n-1}} {\cal{R}}_c\left(\int_0^\infty r^{2n-3} \hat{\phi}(r\cdot) dr\right)(\xi) d\mu(\xi)$$
$$= \int_{S^{2n-1}} \left( \int_{H_\xi} \hat{\phi}(x) dx \right) d\mu(\xi)$$
$$= (2\pi)^{2n-2} \int_{S^{2n-1}} \left( \int_{H_\xi^\bot} \phi(x) dx \right) d\mu(\xi) \ge 0.$$

Now suppose that $\|\cdot\|_K^{-2}$ is a positive definite distribution. By Proposition \ref{posdef},
there exists a finite Borel measure $\mu$ on $S^{2n-1}$ such that for any even test function $\phi$
\begin{equation} \label{eq2}
\int_{\R^{2n}} \|x\|_K^{-2} \phi(x)\ dx = \int_{S^{2n-1}} \left(\int_0^\infty t \hat\phi(t\xi) dt \right) d\mu(\xi).
\end{equation}
Recall that $K$ is $R_\theta$-invariant, so we can assume that $\mu$ is $R_\theta$-invariant and the
latter equality holds only for even $R_\theta$-invariant test functions $\phi.$ For each such test function,
we have by (\ref{eq1}) that
$$\int_0^\infty t \hat{\phi}(t\xi) dt  = \frac1{2\pi}\int_{H_\xi^\bot} \hat{\phi}(x) dx.$$
Using this and Proposition \ref{f-r} and then writing the interior integral in polar coordinates, we get
that  the right-hand side of (\ref{eq2}) is equal to
$$\frac1{2\pi} \int_{S^{2n-1}} \left( \int_{H_\xi^\bot} \hat{\phi}(x) dx \right) d\mu(\xi)=
2\pi \int_{S^{2n-1}} \left( \int_{H_\xi} \phi(x) dx \right) d\mu(\xi)$$
$$= 2\pi  \int_{S^{2n-1}} {\cal{R}}_c\left(\int_0^\infty r^{2n-3} \phi(r\cdot) dr \right)(\xi)\ d\mu(\xi).$$
Writing the left-hand side of (\ref{eq2}) in polar coordinates we get
$$ \int_{S^{2n-1}} \|\theta\|_K^{-2} \left(\int_0^\infty r^{2n-3} \phi(r\theta) dr \right) d\theta$$ 
\begin{equation}\label{eq3}
=2\pi  \int_{S^{2n-1}} {\cal{R}}_c\left(\int_0^\infty r^{2n-3} \phi(r\cdot) dr \right)(\xi)\ d\mu(\xi)
\end{equation}
for any even $R_\theta$-invariant test function $\phi.$ Now let $\phi(x)=u(r)v(\theta)$ for every 
$x\in \R^{2n},$ where $x=r\theta,\ r\in [0,\infty),\ \theta\in S^{2n-1},$ $u\in {\cal{S}}(\R)$ is a non-negative
test function on $\R,$ and $v$ is an arbitrary infinitely differentiable $R_\theta$-invariant function
on $S^{2n-1}.$ Then
$$\int_0^\infty r^{2n-3} \phi(r\theta) dr =  v(\theta) \int_0^\infty r^{2n-3} u(r) dr,$$
so the equality (\ref{eq3}) turns into
$$ \int_{S^{2n-1}} \|\theta\|_K^{-2} v(\theta) d\theta 
=2\pi  \int_{S^{2n-1}} {\cal{R}}_c v(\xi)\ d\mu(\xi).$$
Since infinitely smooth functions are dense in $C_c(S^{2n-1}),$
the latter equality also holds for any function $v\in C_c(S^{2n-1}),$ which means that 
$K$ is a complex intersection body.
\endpf


\section{Characterizations of complex intersection bodies} \label{Char of CIB}

Complex intersection bodies are related to two generalizations of the concept of a real intersection
body. These relations allow to apply to the complex case many results established originally
in the real case.

The concept of a $k$-intersection body was introduced
in \cite{K3,K4}. For an integer $k,\ 1\le k <n$ and star bodies $D,L$ in $\R^n,$
we say that $D$ is the $k$-intersection body of $L$ if for every $(n-k)$-dimensional
subspace $H$ of $\R^n,$
$$|D\cap H^\bot|= |L\cap H|.$$
The class of $k$-intersection bodies was defined in \cite{K4} (see also \cite[Section 4.2]{K1})
as follows. 
\begin{df} Let $1\le k < n.$ We say that an origin symmetric star
body $D$ in $\R^n$ is a $k$-intersection body if there exists 
a finite Borel measure $\mu$ on $S^{n-1}$ so that for every even test function
$\phi$ in $\R^n,$
\begin{equation} \label{defkintgeneral}
\int_{\R^n} \|x\|_D^{-k} \phi(x)\ dx = 
\int_{S^{n-1}} \left(
\int_0^\infty t^{k-1} \hat\phi(t\xi)\ dt\right) d\mu(\xi).
\end{equation}
\end{df}
The class of $k$-intersection bodies is related
to a certain generalization of the Busemann-Petty problem in the same way as
intersection bodies are related to the original problem (see \cite{K1} for details;
this generalization offers a condition that allows to
compare volumes of two bodies in arbitrary dimensions). An equivalent and
probably more geometric way to define $k$-intersection bodies would be to say 
that these bodies are limits in the radial metric of $k$-intersection bodies of star bodies 
(see \cite{EMilman1} or \cite{Ru} for a proof of equivalence of this property to the original 
definition from \cite{K4}).

It was shown in \cite{K4} that an origin symmetric star body $K$ in $\R^n$ is a $k$-intersection body
if and only if the function $\|\cdot\|_K^{-k}$ represents a positive definite distribution. By Theorem \ref{posdef},
\begin{co}\label{Char2} An origin symmetric complex star body $K$ in $\R^{2n}$ is a complex intersection body if
and only if it is a 2-intersection body in $\R^{2n}$ satisfying (\ref{rotation}).
\end{co}

It was proved in \cite[Theorem 3]{KKZ} that every origin symmetric complex convex body $K$ in $\R^{2n}$
is a $(2n-4)$-, $(2n-3)$- and $(2n-2)$-intersection body in $\R^{2n}.$ It follows that
\begin{co}  \label{n<4}Every origin symmetric complex convex body in $\R^6$ and $\R^4$ is a complex intersection body.
\end{co} 

This is no longer true in $\R^{2n},\ n\ge 4$ as shown in \cite[Theorem 4]{KKZ}. The unit balls of complex 
$\ell_q$-balls with $q>2$ are not $k$-intersection bodies for any $1\le k < 2n-4.$

\medbreak

Zhang in \cite{Zh3} introduced another generalization of intersection bodies.
For $1\leq k\leq n-1$, the ($n-k$)-dimensional spherical Radon transform is an operator
${\cal{R}}_{n-k}:C(S^{n-1})\mapsto C(G(n,n-k))$ defined by
$${\cal{R}}_{n-k}(f)(H)=\int_{S^{n-1}\cap H}f(x)dx,\quad H\in G(n,n-k).$$
Here $G(n,n-k)$ is the Grassmanian of $(n-k)$-dimensional subspaces of $\R^n.$
Denote the image of the operator ${\cal{R}}_{n-k}$ by X:
$${\cal{R}}_{n-k}\left(C(S^{n-1})\right)=X\subset C(G(n,n-k)).$$
Let $M^+(X)$ be the space of linear positive continuous functionals on $X$, i.e. for every
$\nu\in M^+(X)$ and non-negative function $f\in X$, we have $\nu(f)\geq0$.

\begin{df}(Zhang \cite{Zh3})
An origin symmetric star body $K$ in $\R^n$ is called a {\it generalized $k$-intersection
body} if there exists a functional $\nu\in M^+(X)$ so that for every $f\in C(S^{n-1})$,
$$\int_{S^{n-1}}\|x\|_K^{-k} f(x)d x= \nu({\cal{R}}_{n-k}{f}).$$
\end{df}

It is easy to see that every complex intersection body in $\R^{2n}$ is a generalized 2-intersection body
in $\R^{2n}.$ If  $K$ is a complex intersection body, then there exists an even $R_\theta$-invariant measure 
$\mu$ on $S^{2n-1}$ such that for every $f\in C(S^{2n-1}),$
$$\int_{S^{2n-1}} \|x\|_K^{-2} f(x) dx = \int_{S^{2n-1}} {\cal{R}}_cf(\xi) d\mu(\xi)$$
$$ =\int_{G(2n,2n-2)} {\cal{R}}_{2n-2} f(H) d\nu(H),$$
where  $\nu$ is a measure on $G(2n,2n-2)$ which is the image of $\mu$ under the mapping $\xi\mapsto H_\xi$ 
from $S^{2n-1}$ to $G(2n,2n-2).$ The measure $\nu$ can be considered as a 
positive continuous functional on $X$ acting by 
$$\nu({\cal{R}}_{n-k}f)= \int_{G(2n,2n-2)} {\cal{R}}_{n-k}f(H) d\nu(H),$$
which implies that $K$ is a generalized 2-intersection body in $\R^{2n}.$

On the other hand,  it was shown in \cite{K4} (see also \cite[Theorem 4.23]{K1}) that every generalized $k$-intersection body in $\R^n$ is a $k$-intersection body. 
So we have shown the following

\begin{pr}
\noindent An origin symmetric complex star body $K$ in $\R^{2n}$ is a complex intersection body if
and only if it is a generalized  $2$-intersection body in $\R^{2n}$ satisfying (\ref{rotation}).
\end{pr}

Let us point out that the latter result is surprising. Combined with Corollary {\ref{Char2}} implies that under the invariance assumption (\ref{rotation}) the class of $2$-intersection bodies is exactly the class of generalized $2$-intersection bodies. Without the invariance assumption (\ref{rotation}) this is no longer true as it follows
from an example of E. Milman {\cite{EMilman2}}. 

\smallskip

Goodey and Weil in {\cite{GW}} proved that any intersection body is the limit (in the radial metric topology) of finite radial sums of ellipsoids. This result has been extended by Grinberg and Zhang \cite{GZ} (see another proof in \cite{EMilman1}) to the case of generalized $k$-intersection bodies  where the radial sum is replaced by the $k$-radial sum. Now we are going to prove a complex version of the result of Goodey and Weil. We do it by adjusting to the complex case the proofs from \cite{EMilman1}
and \cite[Theorem 3.10]{KY}.

We define the complex radial sum $K_{1}+^{c} K_{2}$
of two complex star bodies $K_{1}, K_{2}$ as the complex star body that has radial function 
$$ \rho_{K_{1}+^{c} K_{2}}^{2}  = \rho_{K_{1}}^{2} + \rho_{K_{2}}^{2}  . $$
The latter can be written as
$$\|\cdot\|_{K_{1}+^{c} K_{2}}^{-2} = \|\cdot\|_{K_1}^{-2} + \|\cdot\|_{K_1}^{-2}.$$

\begin{theorem}\label{Thm:approximation}
Let $K$ be an origin symmetric complex star body in $\R^{2n}$. Then $K$ is a complex intersection 
body if and only if $ \|\cdot\|_K^{-2}$ is the limit (in the metric of the space $C_c(S^{2n-1})$) of finite sums
$\|\cdot\|_{E_1}^{-2}+\cdots+ \|\cdot\|_{E_m}^{-2} $, where $E_1$,...,$E_m$ are complex ellipsoids in $\R^{2n}$
(i.e. those ellipsoids in $\R^{2n}$ that are complex convex bodies).
\end{theorem}

In other words,  an origin symmetric  complex star body is an intersection body if and only if it is the
limit (in the radial metric) of complex radial sums of complex ellipsoids.

To prove this result   we need a few lemmas.
For fixed $\xi\in S^{2n-1}$, $a>0$, $b>0$, let $E_{a,b}(\xi)$ be an ellipsoid in $\R^{2n}$ with the norm
$$\|x\|_{E_{a,b}(\xi)}=\left(\frac{( x, \xi )^2+ (x, \xi^\bot)^2}{a^2}+ \frac{|x|_2-( x, \xi )^2- (x,\xi^\bot)^2}{b^2}\right)^{\frac12}, \ x\in \R^{2n}.$$

Clearly, $E_{a,b}(\xi)$ is a complex ellipsoid. In fact, $(x, \xi )^2+ (x, \xi^\bot)^2 = |(x,\xi)_c|^2$
is the modulus squared of the complex scalar product of $x$ and $\xi$ considered
as vectors from $\C^n.$ The latter does not change when $x$ is multiplied by any
complex number of modulus 1, which means that the norm of $E_{a,b}(\xi)$ is
invariant with respect to all rotations $R_\theta.$

Using the formula for the Fourier transform of powers of the Euclidean norm in $\R^{2n}$
(see \cite[p.192]{GS}), we get
$$\left(|x|_2^{-2}\right)^\wedge(\theta)=C(n) |\theta|_2^{ -2n+2},$$
where $C(n) = 2^{2n-3}\pi^n \Gamma(n-1).$
By the connection between linear transformations and the Fourier transform,
\begin{equation}\label{linear transform}
\left( \|Tx\|^{-1}\right)^\wedge(y)=|\det T|^{-1}\left( \|x\|^{-1} \right)^\wedge((T^{*})^{-1} y)
\end{equation}
one can easily compute the following:
\begin{lemma}\label{Lem:FTellipsoid} For all $\theta\in S^{2n-1}$,
$$\left(   \|x\|^{-2}_{E_{a,b}(\xi)}\right)_x^\wedge(\theta) = \frac{C(n)}{a^{2n-4}}
\|\theta\|_{E_{b,a}(\xi)}^{-2n+2}. $$
\end{lemma}
\pf  By (\ref{linear transform}) with $T$ being the composition of a rotation and 
a diagonal operator,
\begin{eqnarray*}
\left(   \|x\|^{-2}_{E_{a,b} (\xi)}\right)_x^\wedge(\theta)  = C(n) {a^2b^{2n-2}}
\|\theta\|_{E_{1/a,1/b}(\xi)}^{-2n+2}=\frac{C(n)}{a^{2n-4}} \|\theta\|_{E_{b,a}(\xi)}^{-2n+2}.
\end{eqnarray*}
 \qed

\begin{lemma}\label{Step1}
Let $K$ be an origin symmetric complex star body, then the function $ \|\xi\|_K^{-2}$ can be  
approximated in the space $C_c(S^{2n-1})$ by functions of the form
\begin{equation}\label{eqn:delta-seq}
f_{a,b}(\xi)=\frac{C(n)}{a^{2n-4}} \int_{S^{n-1}} \|\theta\|_K^{-2} \|\theta\|_{E_{b,a}(\xi)}^{-2n+2}
d\theta,
\end{equation} as $a\to 0$ and $b$ is chosen appropriately.
\end{lemma}

\pf   Using Parseval's formula (\ref{parseval}) and the previous Lemma  we get
$$\frac{C(n)}{  a^{2n-4} } \int_{S^{2n-1}} \|\theta\|_{E_{b,a}(\xi)}^{-2n+2} d\theta$$
$$= \frac{ 1}{  a^{2n-4}  } \int_{S^{2n-1}} \|\theta\|_{E_{b,a}(\xi)}^{-2n+2} \left(|x|_2^{-2}\right)^\wedge(\theta) d\theta$$
 $$= \frac{ 1}{  a^{2n-4}  } \int_{S^{2n-1}} \left(\|x\|_{E_{b,a}(\xi)}^{-2n+2}\right)^\wedge(\theta)
|\theta|_2^{-2}d\theta = \frac{ 1}{C(n)   } \int_{S^{n-1}}   \|x\|_{E_{a,b}(\xi)}^{-2} dx $$
$$= \frac1{C(n)}\int_{S^{2n-1}} \left(\frac{( x, \xi )^2+ (x, \xi^\bot)^2}{a^2}+ \frac{1-( x, \xi )^2- (x,\xi^\bot)^2}{b^2}\right)^{-1} dx.$$

For every fixed $a,$ the latter integral goes to infinity when $b\to \infty,$ and it goes to zero when $b\to 0,$ 
so for every $a>0$ there exists $b=b(a)$ such that
$$\frac{ C(n)}{ a^{2n-4} } \int_{S^{2n-1}} \|\theta\|_{E_{b(a),a}(\xi)}^{-2n+2} d\theta =1.$$
Note that, by rotation invariance, the value of $b(a)$ does not depend on the choice of $\xi.$

Now for every $\xi\in S^{2n-1}$ we have
$$\left|\|\xi\|_K^{-2} - \frac{ C(n)}{ a^{2n-4} } \int_{S^{2n-1}} \|\theta\|_K^{-2}
\|\theta\|_{E_{b(a),a}(\xi)}^{-2n+2} d\theta\right|$$
$$\le \frac{ C(n)}{ a^{2n-4} }\int_{S^{2n-1}}\Big| \|\xi\|_K^{-2} - \|\theta\|_K^{-2}
\Big| \|\theta\|_{E_{b(a),a}(\xi)}^{-2n+2} d\theta$$
$$= \frac{ C(n)}{ a^{2n-4} } \int_{(\theta,\xi)^2+(\theta,\xi^\bot)^2\ge\delta}\Big| \|\xi\|_K^{-2} - \|\theta\|_K^{-2}
\Big| \|\theta\|_{E_{b(a),a}(\xi)}^{-2n+2} d\theta$$
$$+ \frac{ C(n)}{ a^{2n-4} } \int_{(\theta,\xi)^2+(\theta,\xi^\bot)^2<\delta}\Big| \|\xi\|_K^{-2} -  \|\theta\|_K^{-2}
\Big| \|\theta\|_{E_{b(a),a}(\xi)}^{-2n+2} d\theta =I_1+I_2.$$
for every $\delta\in (0,1).$

Since $K$ is a complex star body, the norm of $K$ is constant on vectors of the form 
$u\xi+v\xi^\bot$ with $u^2+v^2=1.$ Vectors of this form are the only solutions on the sphere of the 
equation $(\theta,\xi)^2+(\theta,\xi^\bot)^2=1.$ Using this and 
the uniform continuity of $ \|x\|_K^{-2}$ on the sphere $S^{2n-1}$, for any
given $\epsilon>0$ we can find $\delta\in(0,1)$, close to $1$, so that
$(\theta,\xi)^2+(\theta,\xi^\bot)^2\ge\delta$ implies $\Big| \|\xi\|_K^{-2} - \|\theta\|_K^{-2} \Big|<\epsilon/2$.
Therefore
$$I_1= \frac{ C(n)}{ a^{2n-4} } \int_{(\theta,\xi)^2+(\theta,\xi^\bot)^2\ge\delta}\Big| \|\xi\|_K^{-2} - \|\theta\|_K^{-2}
\Big| \|\theta\|_{E_{b(a),a}(\xi)}^{-2n+2} d\theta$$
$$\le \frac{\epsilon}{2} \left[\frac{ C(n)}{ a^{2n-4} } \int_{(\theta,\xi)^2+(\theta,\xi^\bot)^2\ge\delta}
\|\theta\|_{E_{b(a),a}(\xi)}^{-2n+2} d\theta \right]\le \frac{\epsilon}{2}.$$

Now fix $\delta$ chosen above and estimate the integral $I_2$ as follows
$$I_2=\frac{ C(n)}{ a^{2n-4} } \int_{(\theta,\xi)^2+(\theta,\xi^\bot)^2<\delta}\Big| \|\xi\|_K^{-2} - \|\theta\|_K^{-2}
\Big| \|\theta\|_{E_{b(a),a}(\xi)}^{-2n+2} d\theta$$
$$\le  \frac{ C_1(n,K)}{a^{2n-4}}
\int_{(\theta,\xi)^2+(\theta,\xi^\bot)^2<\delta} \|\theta\|_{E_{b(a),a}(x)}^{-2n+2} d\theta
=\frac{ C_1(n,K)}{a^{2n-4}}\times$$$$\int_{(\theta,\xi)^2+(\theta,\xi^\bot)^2<\delta}  
 \left(\frac{( \theta, \xi )^2+ (\theta, \xi^\bot)^2}{(b(a))^2}+ \frac{1-( \theta, \xi )^2- (\theta,\xi^\bot)^2}{a^2}\right)^{-n+1} d\theta$$
 $$\le a^2 (1-\delta)^{-n+1}C_1(n,K) |S^{2n-1}| .$$
where $$C_1(n,K)={2 C(n) \max_{x\in S^{2n-1}} \|x\|_K^{-2} } .$$ 
Now we can choose $a$ so small that $I_2\le \epsilon/2$. \qed

\bigbreak

\begin{lemma}\label{Step2}
If $\mu$ is a finite measure on $S^{2n-1}$ and $a,b>0$, then the function
$$f(\xi)=\int_{S^{n-1}}   \|\theta \|^{-2}_{E_{a,b}(\xi)} d\mu(\theta)$$ can be  approximated in $C_c(S^{2n-1})$ by the sums of the form
$$\sum_{i=1}^m     \|\xi\|_{E_i}^{-2},$$
where $E_1$,...,$E_m$ are complex ellipsoids.
\end{lemma}

\pf Let $\sigma>0$ be a small number and choose a finite covering of the sphere by spherical
$\sigma$-balls $B_\sigma (\eta_i)=\{\eta\in S^{n-1}: |\eta-\eta_i|<\sigma\}$, $\eta_i\in S^{2n-1}$,
$i=1,\dots, m=m(\delta)$. Define
$$\widetilde{B}_\sigma(\xi_1)=B_\sigma (\xi_1)$$
and
$$\widetilde{B}_\sigma(\xi_i)=B_\sigma (\xi_i)\setminus\bigcup_{j=1}^{i-1}B_\sigma(\xi_j), \quad\mbox{ for } i=2,...,m.$$

Let $1/p_i=\mu(\widetilde{B}_\sigma(\xi_i))$. Clearly, $1/p_1+\cdots+1/p_m=\mu(S^{2n-1})$.

Let $\rho(E_{a,b}(\xi),x)$ be the value of the radial function of the ellipsoid $E_{a,b}(\xi)$ at the point $x$, that is
$$\rho(E_{a,b}(\xi),x)=\|x\|_{E_{a,b}(\xi)}^{-1}.$$

Note that $\rho(E_{a,b}(\xi),x)=\rho(E_{a,b}(x),\xi),$ because both depend only on the modulus of 
the complex scalar product of $x$ and $\xi$, therefore
$$|\rho^2(E_{a,b}(\xi),x)- \rho^2(E_{a,b}(\theta),x)|\le C_{a,b} |\xi - \theta|,$$
with a constant $C_{a,b}$ that depends only on $a$ and $b$.

Then,
$$
\left|\int_{S^{n-1}}  \rho^2(E_{a,b}(\xi),x) d\mu(\xi)-\sum_{i=1}^m\frac{1}{p_i}
\rho^2(E_{a,b}(\xi_i),x)\right|$$
$$ =\left|\sum_{i=1}^m\left(\int_{\widetilde{B}_\sigma(\xi_i)}
\rho^2(E_{a,b}(\xi),x) d\mu(\xi)-\int_{\widetilde{B}_\sigma(\xi_i)} \rho^2(E_{a,b}(\xi_i),x)
d\mu(\xi)\right)\right|$$
$$ \le \sum_{i=1}^m \int_{\widetilde{B}_\sigma(\xi_i)} \left|   {
\rho^2(E_{a,b}(\xi),x)}- {\rho^2(E_{a,b}(\xi_i),x)} \right|d\mu(\xi)$$
$$ \le \sum_{i=1}^m \int_{\widetilde{B}_\sigma(\xi_i)}
  C_{a,b} |\xi-\xi_i|
d\mu(\xi)  \le    C_{a,b}\ \sigma \mu(S^{2n-1}).
$$
Since $\sigma$ is arbitrarily small, the result follows after we define ellipsoids $E_i$ by $$\|x\|_{E_i}^{-2} =
\frac{1}{p_i} \rho^2(E_{a,b}(\xi_i),x). \qed$$
 
\bigbreak

{\sl Proof of Theorem \ref{Thm:approximation}.} The
``if" part immediately follows from Theorem \ref{posdef}, since for any ellipsoid $E$
the distribution $\|\cdot\|_E^{-2}$ is positive definite, as the linear perturbation of the the same 
function for the Euclidean ball.

To prove the converse, suppose that  $K$ is a complex intersection body and $\mu$ is the 
measure on $S^{2n-1}$ corresponding to $K$ by the definition of complex intersection body. 
By Lemma \ref{Step1},  $
\|\xi\|_K^{-2}$ can be uniformly approximated by the integrals of the form
\begin{equation}\label{eq8}
\frac{C(n)}{ a^{2n-4} } \int_{S^{2n-1}} \|\theta\|_K^{-2} \|\theta\|_{E_{b,a}(\xi)}^{-2n+2} d\theta,
\end{equation} as $a\to 0$.

By Lemma \ref{Lem:FTellipsoid} and  Lemma \ref{r-s},
$$\frac{C(n)}{ a^{2n-4} }{\cal{R}}_c(\|\cdot\|_{E_{b,a}(\xi)}^{-2n+2}) = (2\pi)^{2n-1}\|\cdot\|_{E_{a,b}(\xi)}^{-2},$$
and by the definition of complex intersection body, (\ref{eq8}) is equal to
$$(2\pi)^{2n-1} \int_{S^{2n-1}}  \|\theta \|_{E_{a,b}(\xi)}^{-2} d\mu(\theta).$$ Now, by Lemma \ref{Step2},
$\int_{S^{2n-1}}  \|\theta \|_{E_{a,b}(\xi)}^{-2} d\mu(\theta)$ can be uniformly
approximated by sums of the form $\sum_{i=1}^m \|\xi\|_{E_i}^{-1},$
where $E_i$ are complex ellipsoids.
\qed


\section{Stability in the Busemann-Petty problem and hyperplane inequalities} \label{stab-hyper}

Intersection bodies played an important role in the solution
of the Busemann-Petty problem posed in \cite{BP} in 1956.
Suppose that $K$ and $L$ are origin
symmetric convex bodies in $\R^n$ so that, for every $\xi\in S^{n-1},$
$$|K\cap \xi^\bot|\le |L\cap \xi^\bot|.$$
Does it follow that $|K|\le |L|?$
The problem was completely solved at the end of 1990's, and the answer
is affirmative if $n\le 4$ and negative if $n\ge 5.$
The solution appeared as the result of
a sequence of papers \cite{LR}, \cite{Ba}, \cite{Gi}, \cite{Bo4}, 
\cite{L}, \cite{Pa}, \cite{G1}, \cite{G2}, \cite{Z1}, \cite{Z2}, \cite{K5}, \cite{K9}, \cite{Z4},
\cite{GKS} 
(see \cite[Chapter 8]{G} or \cite[Chapter 5]{K1} for details).
One of the main ingredients of the solution was a connection between intersection bodies
and the Busemann-Petty problem established by Lutwak \cite{L}:
if $K$ is an intersection body then the answer to the Busemann-Petty
problem is affirmative for any star body $L.$ On the other hand, if $L$
is a symmetric convex body that is not an intersection body then one can
construct $K$ giving together with $L$ a counterexample.

The complex Busemann-Petty problem can be formulated as follows. Suppose that 
$K,L$ are origin symmetric complex convex bodies in $\R^{2n}$ and,  for every
$\xi \in S^{2n-1},$ we have $|K\cap H_\xi|\le |L\cap H_\xi|.$ Does it follow that $|K|\le |L|?$
As proved in \cite{KKZ}, the answer is affirmative if $n\le 3,$ and it is negative if $n\ge 4.$
The proof is based on a connection with intersection bodies, similar to Lutwak's connection
in the real case (see \cite[Theorem 2]{KKZ}):  

\noindent (i) If $K$ is a complex intersection body in $\R^{2n}$ 
and $L$ is any origin symmetric complex star body in $\R^{2n},$ then the answer to the question 
of the complex Busemann-Petty problem is affirmative;

\noindent (ii) if there exists an origin symmetric complex convex body in $\R^{2n}$ that is not a complex intersection
body , then one can construct a counterexample to the complex Busemann-Petty problem.

These connections were formulated in \cite{KKZ} in terms of positive definite distributions, so one has to use 
Theorem \ref{posdef} to get the statements in terms of convex intersection bodies. 

Zvavitch \cite{Zv} found a generalization of the Busemann-Petty problem
to arbitrary measures, namely, one can replace volume by any measure $\gamma$ with 
even continuous density in $\R^n.$ In particular, if $n\le 4,$ then for any 
origin symmetric convex bodies $K$ and $L$ in $\R^n$ the inequalities
$$\gamma(K\cap \xi^\bot) \le \gamma(L\cap \xi^\bot), \qquad \forall \xi\in S^{n-1}$$
imply
$$\gamma(K)\le \gamma(L). $$
Zvavitch also proved that this is generally not true if $n\ge 5,$ namely, for any $\gamma$ with strictly positive
even continuous density there exist $K$ and $L$ providing a counterexample. 
In \cite{Zy} the result of Zvavitch was extended to complex convex bodies.

In this section we are going to prove stability in the affirmative part of the result from \cite{Zy}; 
see Theorem \ref{stab-complex} below. Note that stability in the original Busemann-Petty problem 
was established in \cite{K7}, and for the complex Busemann-Petty problem it was done in \cite{K6}. 
Stability in Zvavitch's result was proved in \cite{K8}, and in \cite{KM} the result of \cite{K8} was extended 
to sections of lower dimensions in place of hyperplane sections.

Let $f$ be an even continuous non-negative function on $\R^{2n},$ and denote by $\gamma$ the measure on $\R^{2n}$ with density $f$
so that for every closed bounded set $B\subset \R^n$
$$\gamma(B)=\int_B f(x)\ dx.$$ 
Since we apply $\gamma$ only to complex star bodies, we can assume without loss of generality
that the measure $\gamma$ and the function $f$ are $R_\theta$-invariant.

We need a polar formula for the measure of a complex star body $K:$
\begin{equation}\label{meas-polar}
\gamma(K) = \int_K f(x)\ dx = \int_{S^{n-1}} \left(\int_0^{\|\theta\|_K^{-1}} r^{n-1}f(r\theta)\ dr \right) d\theta.
\end{equation}

For every $\xi\in S^{2n-1},$
$$ \gamma(K\cap H_\xi) = \int_{K\cap H_\xi} f(x) dx$$
$$= \int_{S^{2n-1}\cap H_\xi} \left(\int_0^{\|\theta\|_K^{-1}}r^{2n-3} f(r\theta) dr \right)d\theta$$
\begin{equation} \label{meas-sect-polar}
= {\cal{R}}_c \left(\int_0^{\|\cdot\|_K^{-1}} r^{2n-3} f(r\cdot)\ dr \right)(\xi),
\end{equation}

We need the following elementary lemma, which was also used by Zvavitch in \cite{Zv}.
\begin{lemma}\label{zvavitch} Let $a,b \in [0,\infty),\ n\in \N, n\ge 2$, and 
let $g$ be a non-negative integrable function on $[0,\max\{a,b\}].$ Then
$$\int_0^a r^{2n-1} g(r)\ dt-a^2\int_0^ar^{2n-3} g(r)\ dr$$
$$\le \int_0^b r^{2n-1} g(r)\ dr - a^2\int_0^br^{2n-3} g(r)\ dr.$$
\end{lemma}

Denote by  $$d_n= \frac{|B_2^{2n}|^{\frac{n-1}n}}{|B_2^{2n-2}|},$$ where $B_2^n$ stands for the unit Euclidean ball in $\R^n,$
and $$|B_2^n|=\frac{\pi^{n/2}}{\Gamma(1+\frac n2)}.$$
Note that $d_n<1$ for every $n\in \N;$ this easily follows from the log-convexity of the $\Gamma$-function.

It is well-known (see for example \cite[p.32]{K1}) that the surface area of the sphere $S^{n-1}$ in $\R^n$ is equal to 
\begin{equation} \label{ball}
|S^{n-1}| = n |B_2^n|.
\end{equation}

The following stability result extends \cite[Theorem 2]{K8}  to the complex case. The proof is similar to that of \cite[Theorem 2]{KM}.
\begin{theorem}\label{stab-complex} Let $K$ and $L$ be origin symmetric complex star bodies in $\R^{2n},$
let $\e>0$ and let $\gamma$ be a measure on $\R^{2n}$with even continuous non-negative density $f.$ 
Suppose that $K$ is a complex intersection body, and that for every $\xi\in S^{2n-1}$
\begin{equation}\label{stab-ineq}
\gamma(K\cap H_\xi) \le \gamma(L\cap H_\xi) + \e.
\end{equation}
Then 
$$\gamma(K)\le \gamma(L) + \frac n{n-1} d_n \e |K|^{\frac 1n}.$$
\end{theorem}
\pf  By (\ref{meas-sect-polar}), the condition (\ref{stab-ineq}) can be written as 
$${\cal{R}}_c \left(\int_0^{\|\cdot\|_K^{-1}} r^{2n-3} f(r\cdot)\ dr \right)(\xi) $$
$$\le  {\cal{R}}_c \left(\int_0^{\|\cdot\|_L^{-1}} r^{2n-3} f(r\cdot)\ dr \right)(\xi) + \e, \qquad \forall \xi\in S^{2n-1}.$$
Integrating the latter inequality with respect to the measure $\mu$ on $S^{2n-1},$ corresponding to the body $K$
by Definition \ref{int}, and then using the equality of Definition \ref{int}, we get
\begin{equation}\label{eq6}
\int_{S^{n-1}} \|\theta\|_K^{-2} \left(\int_0^{\|\theta\|_K^{-1}} r^{2n-3} f(r\theta)\ dr\right) d\theta 
\end{equation}
$$\le \int_{S^{2n-1}} \|\theta\|_L^{-2} \left(\int_0^{\|\theta\|_L^{-1}} r^{2n-3} f(r\theta)\ dr\right) d\theta
    + \e \int_{S^{2n-1}} d\mu(\xi).$$
Applying Lemma \ref{zvavitch} with $a=\|\theta\|_K^{-1}$, $b=\|\theta\|_L^{-1}$
    and $g(r)=f(r\theta)$ and then integrating over the sphere, we get
  
$$ \int_{S^{2n-1}} \left(\int_0^{\|\theta\|_K^{-1}} r^{2n-1} f(r\theta)\ dr \right) d\theta$$ 
$$- \int_{S^{2n-1}} \|\theta\|_K^{-k} \left(\int_0^{\|\theta\|_K^{-1}} r^{2n-3} f(r\theta)\ dr \right) d\theta$$
$$\le \int_{S^{2n-1}} \left(\int_0^{\|\theta\|_L^{-1}} r^{2n-1} f(r\theta)\ dr \right) d\theta$$
\begin{equation} \label{eq7}
- \int_{S^{2n-1}} \|\theta\|_K^{-k} \left(\int_0^{\|\theta\|_L^{-1}} r^{2n-3} f(r\theta)\ dr \right) d\theta.
\end{equation}
Adding (\ref{eq6}) and (\ref{eq7}) and using (\ref{meas-polar}) we get
$$\gamma(K)\leq\gamma(L)+\e \int_{S^{2n-1}} d\mu(\xi).$$
Since ${\cal{R}}_c 1 = |S^{2n-3}| 1,$ where $1(\xi)\equiv 1,$ we again apply Definition \ref{int},
H\"older's inequality, the polar formula for volume and (\ref{ball}):
$$\int_{S^{2n-1}} d\mu(\xi) = \frac 1{|S^{2n-3}|} \int_{S^{2n-1}} {\cal{R}}_c1(\xi) d\mu(\xi)$$
$$= \frac 1{|S^{2n-3}|} \int_{S^{2n-1}}  \|x\|_K^{-2} dx \le  \frac 1{|S^{2n-3}|} \left(\int_{S^{2n-1}} \|x\|_K^{-2n}dx \right)^{\frac 1n}
|S^{2n-1}|^{\frac{n-1}n}$$
$$= \frac{(2n)^{\frac 1n}|S^{2n-1}|^{\frac{n-1}n}}{|S^{2n-3}|} |K|^{\frac 1n} = \frac n{n-1} d_n |K|^{\frac1n}.$$
\endpf

Interchanging $K$ and $L$ in Theorem \ref{stab-complex}, we get a complex version of \cite[Corollary 1]{K8}.
\begin{co} If $K$ and $L$ are complex intersection bodies in $\R^{2n},$ then
$$|\gamma(K) - \gamma(L)|$$
$$\le \frac n{n-1} d_n \max_{\xi\in S^{2n-1}} |\gamma(K\cap H_\xi)-\gamma(L\cap H_\xi)|
\max \left\{|K|^{\frac1n}, |L|^{\frac 1n} \right\}.$$
\end{co}

Putting $L=\emptyset$ in the latter inequality, we extend to the complex case the hyperplane inequality
for real intersection bodies from \cite[Theorem 1]{K8}.

\begin{theorem} \label{hyper-ineq} If $K$ is a complex intersection body in $\R^{2n}$, and $\gamma$ 
is an arbitrary measure on $\R^{2n}$ with even continuous density, then
$$\gamma(K) \le \frac n{n-1} d_n \max_{\xi\in S^{2n-1}} \gamma(K\cap H_\xi)\ |K|^{\frac1n}.$$
By Corollary \ref{n<4}, this inequality holds for any origin symmetric complex convex body $K$ in $\R^4$ or $\R^6.$
\end{theorem}

The constant in Theorem \ref{hyper-ineq} is optimal, as can be easily seen from the same example as in \cite{K8}. Let $K=B_2^n$
and, for every $j\in N,$  let $f_j$ be a non-negative continuous function on $[0,1]$ supported in $(1-\frac 1j,1)$
and such that $\int_0^1 f_j(t) dt =1.$ Let $\gamma_j$ be the measure on $\R^{2n}$ with density $f_j(|x|_2),$
where $|x|_2$ is the Euclidean norm in $\R^{2n}.$  Then a simple computation shows that 
$$\lim_{j\to \infty} \frac{\gamma_j(B_2^{2n})}{\max_{\xi\in S^{2n-1}} \gamma_j(B_2^{2n}\cap H_\xi)\ |B_2^{2n}|^{1/n}}
= \frac n{n-1}d_n.$$

Note that in the case of volume (when the density $f\equiv 1$) the inequality of Theorem \ref{hyper-ineq} follows
from \cite[Corollary 1]{K6} and the constant is just $d_n$ without the term $n/(n-1).$ One has to follow the proof of 
Theorem 1 from \cite{K6} to restore the constant $d_n$ which is estimated by 1 everywhere in \cite{K6}. 

The result of Theorem \ref{hyper-ineq} is related to the famous  hyperplane problem asking whether there exists 
an absolute constant $C$ so that for any origin symmetric convex body $K$ in $\R^n$
\begin{equation} \label{hyper}
\vol_n(K)^{\frac {n-1}n} \le C \max_{\xi \in S^{n-1}} \vol_{n-1}(K\cap \xi^\bot),
\end{equation}
where  $\xi^\bot$ is the central hyperplane in $\R^n$ perpendicular to $\xi.$
The problem is still open, with the best-to-date estimate $C\sim n^{1/4}$ established
by Klartag \cite{Kl}, who slightly improved the previous estimate of Bourgain \cite{Bo3}.


 \section{Complex  intersection bodies of convex complex bodies} \label{bus}

 In this section we extend two classical results about intersection bodies of convex bodies to the complex setting. 
 The well-known result of Busemann \cite{Bu1} is that the intersection body of a symmetric convex
 body is also symmetric convex. We prove a complex version of this result.
 
  \begin{theorem}\label{CBusemann}
Let $K$ be an origin symmetric convex body in $\C^n$ and $ I_{c}(K)$ the complex intersection body of $K.$ Then $ I_{c}(K)$ is also an origin symmetric  convex body in $\C^{n}$.
\end{theorem}

 Before we prove Theorem \ref{CBusemann} we need some preparations.
 We write $O_{2}$ for the set or all rotations in $\R^{2}$ as described in \eqref{rotation} and   $H_{0}$ for the (real) $2$-dimensional subspace of $\R^{2n}$, spanned by the standard unit vectors $e_{1}, e_{2}$.  If $V_{H_{0}}$ is any orthogonal transformation in $H_{0}$ ($V_{H_{0}} \in O_{2}$), we define $U_{V_{H_0}}\in O_{2n}$ to be the orthogonal transformation in $\R^{2n}$ that has (as a matrix) $n$ copies of $V_{H_{0}}$ on its diagonal. Then \eqref{rotation} implies that $ U_{V_{H_0}} K = K$. We will use this property in the following form:
\begin{equation}\label{characterization} {\bf 1}_{K} ( U_{V_{H_0}}x)= {\bf 1}_{K} (x) , \ \forall x\in \R^{2n} , \ \forall  \ V_{H_0} \in O_2. \end{equation}

\noindent Assume that $n\geq 3$. Let $u_{1}, u_{2} \in \C^n, \  |u_1|_{2}=|u_2|_{2}=1,$ with $H_{u_i}^{\perp}= {\rm span} \{u_i, \ u_i^{\perp}\},$ as introduced earlier, and $\theta_i\in S_{H_{u_i}^{\perp}}, \  i=1, 2.$ We define $u_{3}:= \frac{ u_{1}+u_{2}}{ |u_{1}+ u_{2}|_{2}} ,$ with $H^{\perp}_{u_3}$ and $\theta_{3}:= \frac{ \theta_{1}+ \theta_{2}}{ |\theta_{1}+ \theta_{2}|_{2}}\in S_{H^{\perp}_{u_3}}$ such that $|\theta_{1}+ \theta_{2}|_{2}=|u_{1}+ u_{2}|_{2}$. We can assume that $H_{u_{1}}^{\perp}\cap H_{u_{2}}^{\perp}=\{0\}$. 

\noindent Now, let $r_{1}, r_{2}>0$. We define $r_{3}, t$ (as functions of $r_{1}, r_{2}$) such that 
\begin{equation}\label{r,t} t:= \frac{r_{1}}{ r_{1}+ r_{2}} , \ r_{3} :=  \frac{1}{ \frac{1}{r_{1}} + \frac{1}{r_{2}} }| u_{1}+ u_{2}|_{2} .\end{equation}

If $S:=(H_{u_1}^{\perp} + H_{u_2}^{\perp})^{\perp},$ we write $ E_{i}:= {\rm span}\{ H_{u_{i}}^{\perp}, S\}, \ i=1,2,3$.  We define the functions $g_{i}:H_{u_i}^{\perp} \rightarrow \R$, $h_{i}:[0,\infty)\rightarrow [0,\infty), \ i=1,2,3$ to be 
\begin{equation}\label{g,h} g_{i}(x):= \int_{S+x} {\bf 1}_{K}(y) d y = |K \cap ( S+ x)| , \ h_{i}(r):= g_{i}(r\theta_{i} ).\end{equation}

 In the following we exploit the fact that $K$ satisfies \eqref{characterization}.
 
 \begin{lemma}\label{l2.1}
For $i=1,2,3$, we have that 
\begin{equation}\label{eqKEi}| K\cap E_{i}| = 2\pi \int_{0}^{\infty} r h_{i}(r) d r .\end{equation}
\end{lemma}

  \pf
  First we will show that the functions $g_{i}$ are rotation invariant. Let $S_{0}:= {\rm span}\{ e_{5}, \cdots, e_{2n}\}$ and $\theta_{i}^{(1)}, \theta_{i}^{(2)} \in S_{H^{\perp}_{U_i}} \ i=1,2,3.$ There exists $ U_{i}\in O_{2n}$ such that $ U_{i}H^{\perp}_{u_i} = H_{0}$ and $ U_{i}S= S_{0},$ and $\phi_{i}^{(1)}, \phi_{i}^{(2)} \in S_{H_0},$ such that 
  $U^{\top}_i (\phi_{i}^{(1)})=\theta_{i}^{(1)}$ and $U^{\top}_i (\phi_{i}^{(2)})=\theta_{i}^{(2)}, \ i=1,2,3.$ Moreover, there exists $V\in O_2$ such that $V(\phi_i^{(2)})=\phi_i^{(1)}.$ Let $V_0\in O_{2n}$ be the diagonal operator with $V$ on its diagonal entries. Then, it is clear that  $V_0(\phi_{i}^{(2)})=\phi_{i}^{(1)}$ and $V_0S_0=S_0.$ Then, by \eqref{characterization}, we have that for every $r>0,$
  
  $$ g_{i}(r\theta_{i}^{(1)}) =  \int_{S+r\theta_{i}^{(1)}} {\bf 1}_{K}(y) d y = \int_{U_{i}^{\top} S_{0} + r U_{i}^{\top} (\phi_{i}^{(1)}) } {\bf 1}_{K} (y ) d y $$
$$ =  \int_{U_{i}^{\top} (S_{0} + r   \phi_{i}^{(1)} )} {\bf 1}_{K} (y ) d y =  \int_{S_{0} + r\phi_{i}^{(1)} }  {\bf 1}_{K} ( U_{i} y ) d y $$
$$= \int_{V_{0} S_{0} + rV_{0}(\phi_{i}^{(2)}) }  {\bf 1}_{K} ( U_{i} y ) d y =  \int_{V_{0} (S_{0} + r\phi_{i}^{(2)}) }  {\bf 1}_{K} ( U_{i} y ) d y  $$
$$ = \int_{  S_{0} + r \phi_{i}^{(2)} }  {\bf 1}_{K} ( V^{\top}_{0} U_{i} y ) d y =    \int_{  S_{0} + r\phi_{i}^{(2)} }  {\bf 1}_{K} (   U_{i} y ) d y $$
$$  = \int_{   U_{i}^{\top} S_{0} + r  U_{i}^{\top}( \phi_{i}^{(2)}) }  {\bf 1}_{K} (    y ) d y =  \int_{  S + r   \theta_{i}^{(2)} }  {\bf 1}_{K} (    y ) d y = g_{i}(r\theta_{i}^{(2)})  . $$
So, by Fubini's theorem we have that for $i=1,2,3$,
$$ |K\cap E_{i}|= \int_{H_{u_{i}}^{\perp}} \int_{S+ x} {\bf 1}_{K}(y) d y dx = \int_{H_{u_{i}}^{\perp}} g_{i}(x) d x  $$
$$ =  \int_{S_{H_{u_{i}}^{\perp}}} \int_{0}^{\infty} r g_{i}(r\theta) d r d\theta =  2\pi  \int_{0}^{\infty} r g_{i}(r\theta_{i}) d r  = 2\pi  \int_{0}^{\infty} r h_{i}(r ) d r.   $$
This finishes the proof. $\hfill\Box $ 
  
\medbreak  
\noindent The convexity of $K$ is exploited in the following 

\begin{lemma}{\label{l2.2}}
\noindent With the above notation we have that 
\begin{equation}\label{2.6} h_{3}(r_{3}) \geq h_{1}^{(1-t)}(r_{1})h_{2}^{t}(r_{2}) . \end{equation}
\end{lemma}

  \pf
\noindent Note that $ h_{i}(r):= | K \cap ( S + r\theta_{i})|$, for $i=1,2,3$. Also observe that the sets $ K\cap (S + r_{i} \theta_{i})$ all lie in the same hyperplane, and by convexity, 
$$ (1-t) (K\cap (S+ r_{1}\theta_{1}) ) + t  (K\cap (S+ r_{2}\theta_{2}) )\subseteq  (K\cap (S+ r_{3}\theta_{3}) ) . $$
The result follows from the Brunn-Minkowski inequality. $\hfill\Box $ 

\smallskip

\noindent Finally we need the following result of K. Ball \cite{Ball1}. It can be seen as an extension of the inequality of Busemann (see also \cite{BFranz}). In \cite{Ball1} this proposition has been proved but has not been stated in this form. In this form it can be found in \cite{Kl}. 

\begin{pr}{\label{pr2.3}}
\noindent Let $r_{1}, r_{2}>0$. Define $t, r_{3}$ as follows:
\begin{equation}\label{2.7}  t:= \frac{r_{1}}{ r_{1}+ r_{2}} , \ r_{3} :=  \frac{2}{ \frac{1}{r_{1}} + \frac{1}{r_{2}} } .  \end{equation}
Assume that  $h_{1}, h_{2}, h_{3}:[0,\infty)\rightarrow [0,\infty)$ such that 
\begin{equation}\label{2.8} h_{3}(r_{3}) \geq h_{1}^{(1-t)}(r_{1})h_{2}^{t}(r_{2}) , \ \forall r_{1}, r_{2}>0 . \end{equation}
Let $p\geq 1$ and denote 
$$ A:= \left(\int_{0}^{\infty} r^{p-1} h_{1}(r) d r\right)^{\frac{1}{p}} ,$$
$$ B:=\left( \int_{0}^{\infty} r^{p-1} h_{2}(r) d r\right)^{\frac{1}{p}}  ,$$
$$ C:= \left(\int_{0}^{\infty} r^{p-1} h_{3}(r) d r\right)^{\frac{1}{p}}  . $$
Then,
\begin{equation}\label{2.10} C\geq \frac{2}{ \frac{1}{A}+ \frac{1}{B}} .  \end{equation}
\end{pr}

\smallskip

\noindent We rewrite the previous proposition in a form that fits our setting:

\begin{co}{\label{co2.4}}
\noindent Let $r_{1}, r_{2}>0$ and let $\alpha>0$. Define $t, r_{3}$ as follows:
\begin{equation}\label{2.11}   t:= \frac{r_{1}}{ r_{1}+ r_{2}} , \ r_{3} :=  \frac{\alpha}{ \frac{1}{r_{1}} + \frac{1}{r_{2}} } .  \end{equation}
Assume that  $h_{1}, h_{2}, h_{3}:[0,\infty)\rightarrow [0,\infty)$ such that 
\begin{equation}\label{2.12}  h_{3}(r_{3}) \geq h_{1}^{(1-t)}(r_{1})h_{2}^{t}(r_{2}) , \ \forall r_{1}, r_{2}>0  . \end{equation}
Let $p\geq 1$ and denote 
$$ A:= \left(\int_{0}^{\infty} r^{p-1} h_{1}(r) d r\right)^{\frac{1}{p}}  , $$
$$ B:= \left(\int_{0}^{\infty} r^{p-1} h_{2}(r) d r\right)^{\frac{1}{p}}  , $$
$$ C:= \left(\int_{0}^{\infty} r^{p-1} h_{3}(r) d r\right)^{\frac{1}{p}}  .$$
Then,
\begin{equation}\label{2.14}  C\geq \frac{\alpha}{ \frac{1}{A}+ \frac{1}{B}} . \end{equation}
\end{co}

  \pf
\noindent Let $H_{3}:[0,\infty)\rightarrow [0,\infty)$ be such that $ H_{3}( r) = h_{3}(\frac{\alpha}{2} r),$ and let $ C^{\prime}:= \int_{0}^{\infty} r^{p-1} H_{3}(r) d r$. Let $r_{3}^{\prime}:= \frac{2}{\alpha} r_{3}$. Then $r_{1}, r_{2}, r_{3}^{\prime}$ satisfy ({\ref{2.7}}) and 
$$ H_{3}(r_{3}^{\prime}) = h_{3}(r_{3}) \geq h_{1}^{(1-t)}(r_{1})h_{2}^{t}(r_{2}) , $$
so $h_{1}, h_{2}, H_{3}$ satisfy also ({\ref{2.8}}). So by Proposition \ref{pr2.3} we have that 
$$   \frac{2}{ \frac{1}{A}+ \frac{1}{B}} \leq C^{\prime} = \left(\int_{0}^{\infty} r^{p-1} h_{3}\left(\frac{\alpha}{2}r\right) d r\right)^{\frac{1}{p}}  =  \frac{2}{\alpha} \left( \int_{0}^{\infty} s^{p-1} h_{3}(s) d s\right)^{\frac{1}{p}}  $$

$$ {\rm or} \ \ \ \ C\geq   \frac{\alpha}{ \frac{1}{A}+ \frac{1}{B}} . $$
This completes the proof. $\hfill\Box $ 

\smallskip

\begin{pr}{\label{pr2.5}}
\noindent In the notation introduced above, if $\alpha := | \theta_{1} + \theta_{2}|_{2} $ then    
\begin{equation}\label{2.15} \frac{\alpha}{|K\cap E_{3}|^{\frac{1}{2}}} \leq \frac{1}{|K\cap E_{1}|^{\frac{1}{2}}} + \frac{1}{|K\cap E_{2}|^{\frac{1}{2}}} .\end{equation}
\end{pr}

  \pf
\noindent By Lemma \ref{l2.2} we have that $h_{1}, h_{2}, h_{3}, r_{1}, r_{2}, r_{3},t$ satisfy (\ref{2.11}) and (\ref{2.12}). By Corollary \ref{co2.4}, applied to $p=2,$ we have that for $A,B,C$ as in Corollary \ref{co2.4},
\begin{equation}\label{2.16}  \frac{\alpha}{C}\leq \frac{1}{A}+ \frac{1}{B} . \end{equation}
Note that by Lemma {\ref{l2.1}}, 
$$  A^{2}:= \int_{0}^{\infty} r h_{1}(r) d r = \frac{|K\cap E_{1}|}{2\pi}, \ B^{2}:= \int_{0}^{\infty} r h_{2}(r) d r= \frac{|K\cap E_{2}|}{2\pi} \ {\rm and} $$
\begin{equation}\label{2.17}  C^{2}:= \int_{0}^{\infty} r^{p-1} h_{3}(r) d r= \frac{|K\cap E_{3}|}{2\pi} . \end{equation}
By (\ref{2.16}) and (\ref{2.17}), we complete the proof. $\hfill\Box $ 

\smallskip

\begin{co}{\label{co2.6}}
\noindent Let $K$ be a symmetric complex convex body in $\R^n$. Let $H$ be an $(n-2)$-dimensional subspace of $\C^{n}$. Let $u\in H^{\perp}$ complex unit vector and let $H_{u}:= {\rm span}\{ H, u\}$ and $ r(u):= | K \cap H_{u}|^{\frac{1}{2}}$.
Then $r:H^{\perp} \cap S^{2n-1}\rightarrow (0,\infty)$ is the boundary of a complex convex body in $H^{\perp}$.   
\end{co}

  \pf
\noindent  In order to show that the curve $r$ is the boundary of a convex body in $\C^{n}$ it is enough to show that $r^{-1}$ is the restriction of a norm to $H^{\perp}$. So if $u_{1}, u_{2}$ are two non-parallel unit vectors in $H^{\perp},$ let $u_{3}:= \frac{u_{1}+u_{2}}{ |u_{1}+ u_{2}|_{2}}=: \frac{1}{\alpha} (u_{1}+ u_{2})$. It is enough to show that
$$ \frac{\alpha}{r(u_{3})} \leq \frac{1}{r(u_{1})} +\frac{1}{r(u_{2})} . $$
In the notation of this section we have that $ H_{u_{i}}=E_{i}$ and $ r(u_{i}) = |K \cap E_{i}|^{\frac{1}{2}}$. The result follows from Proposition {\ref{pr2.5}}.  $\hfill\Box $ 

\medskip

{\sl Proof of Theorem {\ref{CBusemann}}:} In the case where $n=2$ the body $I_{c}(K)$ is simply a rotation of $K$, so the result is obvious. Let $n\geq 3$. Then Corollary \ref{co2.6} implies that $I_{c}(K)\cap H^{\perp}$ is convex for every $(n-2)$-dimensional subspace $H$ of $\C^{n}$. This implies that $ I_{c}(K)$ is convex. The symmetry of $I_{c}(K)$ is obvious from the definition. Finally, it is not difficult to see that $I_{c}(K)$ satisfies (\ref{rotation}). This implies that $I_{c}(K)$ is a complex convex body.  $\hfill\Box $

\medskip

In the case where $K\subset \R^n$ is convex, by results of  Hensley \cite{H} and Borell \cite{Bor}, one has that the intersection body of $K$, $I(K)$, is isomorphic to an ellipsoid; i.e.  $ d_{BM}(I(K), B_{2}^{n})\leq c$ where $d_{BM}$ stands for the Banach-Mazur distance and $c>0$ is a universal constant.  Recall that  the Banach-Mazur distance of two symmetric convex bodies $K_{1}, K_2$ in $\R^n$ as 
$$ d_{BM}(K_{1}, K_{2}) := \inf_{T\in GL_{n}} \inf\{ a>0 : K_{1}\subseteq T K_{2} \subseteq a K_{1}\} . $$ 
One can show that the same result holds also in the complex case by using a result of K. Ball {\cite{Ball1}}. However, we can immediately deduce the ``complex Hensley" theorem by using a more general fact (where the result of K. Ball has been used) proved in \cite{KPZ}: 

\begin{pr}
\noindent Let $K$ be an origin symmetric convex body in $\R^n$ and assume that the $2$-intersection body, $ I_{2}(K)$, exists and it is convex. Then 
$$ d_{BM} (I_{2}(K), B_{2}^{n})\leq c , $$
where $c>0$ is an absolute constant. 
\end{pr}

Combining the above result with Proposition {\ref{Char2}} and Theorem {\ref{CBusemann}}, we immediately get the following

\begin{theorem}
\noindent Let $K$ be an origin symmetric convex body in $\C^{n}$. Then 
$$ d_{BM}( I_{c}(K), B_{2}^{n}(\C)) \leq c, $$
where $c>0$ is an absolute constant and $ B_{2}^{n}(\C):= B_{2}^{2n}$ is the Euclidean ball in $\C^{n}$.  
\end{theorem}

\bigskip

{\sl Acknowledgments:} The first named author wishes to thank the US National Science Foundation for support through grant
DMS -1001234. The second named author wishes to thank the A. Sloan Foundation and US National Science Foundation for support through grant DMS-0906051 . Part of this work was carried out when the third named author was visiting the Mathematics Department of Texas A$\&$M University, which she thanks for hospitality.


\footnotesize
\bibliographystyle{amsplain}

\end{document}